\documentclass[12pt]{article}
\usepackage{graphicx}
\usepackage{amsmath}
\usepackage{amsthm}
\usepackage{latexsym,bm}
\usepackage{amssymb}
\usepackage{indentfirst}
\newtheorem{thm}{Theorem}[section]
\newtheorem{cor}[thm]{Corollary}

\newtheorem{prop}[thm]{Proposition}
\newtheorem{rem}[thm]{Remark}
\newtheorem{exm}{Example}

\begin{document}
\title{{Bernstein Type Theorems For Minimal Lagrangian Graphs of Quaternion Euclidean space}
\footnotetext{{ Mathematics Classification Primary(2000)}:{
53A10,53A07,53C38}\\
\hspace*{7mm}{Supported by the Zhongdian grant of NSFC}.\\
\hspace*{7mm}{Keywords: Bernstein type theorem, Quaternion
Euclidean space, minimal Lagrangian graphs.\\
}}}
\author{Yuxin Dong,  Yingbo Han and  Qingchun Ji}

\date{}
\maketitle
\begin{abstract}
In this paper, we prove some Bernstein type results for
$n$-dimensional minimal Lagrangian graphs in quaternion Euclidean
space $H^n\cong R^{4n}$. In particular, we also get a new
Bernstein Theorem for special Lagrangian graphs in $C^n$.
\end{abstract}

\section{Introduction}
The celebrated theorem of Bernstein says that the only entire
minimal graphs in Euclidean 3-space are planes. This result has
been generalized to $R^n$, for $n\leq 7$ and general dimension
under various growth condition, see \cite{han7} and the reference
therein for codimension one case. For higher codimension, the
situation becomes more complicated. Due to the counterexample of
Lawson-Osserman \cite{han0}, the higher codimension Bernstein type
result is not expected to be true in the most generality. Hence we
have to consider the additional suitable conditions to establish a
Bernstein type result for higher codimension.

In recent years, remarkable progress has been made by \cite{han1},
\cite{han2}, \cite{han3}, \cite{han5} and \cite{han6} in Bernstein
type problems of minimal submanifolds with higher codimension and
special Lagrangian submanifolds. The key idea in these papers is
to find a suitable subharmonic function, whose vanishing implies
the minimal graph is totally geodesic. Let $M$ be a minimal
submanifold of $R^{n+m}$ that can be represented as the graph of a
smooth map $f:R^n\longrightarrow R^m$. The function is given by
$$*\Omega=\frac{1}{\sqrt{det(I+(df)^tdf)}}$$
Jost-Xin \cite{han1} established a Bernstein result for $M$ under
the condition $*\Omega\geq K>\frac{1}{2}$, which improves the
previous results in \cite{han10} and \cite{han9}. Wang in
\cite{han5} derived a nice Bochner type formula for the function
$ln(*\Omega)^{-1}$. Under the so-called area-decreasing condition,
he obtained a Bernstein result for higher codimension case too.

Due to string theory, special Lagrangian submanifolds received
much attention in recent years. Some authors also tried to
establish Bernstein type results for special Lagrangian
submanifolds (see \cite{han2}, \cite{han3} and \cite{han6}). It is
known that a special Lagrangian graph may be represented by a
gradient of a smooth function, i.e., $f=\nabla u$ for a smooth
function $u:R^n \longrightarrow R$. The function $u$ is called the
potential function. Tsui and Wang in \cite{han3} obtained
Bernstein results for special Lagrangian graphs by applying the
same Bochner formula. We should point out that the same formula
for special Lagrangian graphs was also derived by Yuan in
\cite{han6} from a different point of view. An important technique
used by Yuan is the so-called Lewy transformation which allows him
to prove: Any special Lagrangian graph given by a convex potential
function must be an affine plane. Actually, Yuan obtained a little
bit stronger Bernstein result for special Lagrangian graph under
the condition $Hess(u)\geq-\epsilon(n)I$, where $\epsilon(n)$ is a
small dimensional constant.

In this paper, we will investigate a real minimal Lagrangian graph
$\Sigma^n$ in Quaternion space $H^n\cong R^{4n}$ which is given by
three potential functions $u_s:R^n\longrightarrow R$, $s=1,2,3$ as
follows:
$$\Sigma=\{(x,\nabla u_1,\nabla u_2, \nabla u_3):x\in R^n\}$$
The Lagrangian condition forces the three matrices $Hess(u_s)$ to
be commutative with each other. As a result, we may choose a
particular quaternionic frame corresponding to singular value
decomposition of $\nabla u_s$ $(s=1,2,3)$ at each point. A useful
formula for the minimal Lagrangian graph is derived by applying
Wang's Bochner formula to the quaternionic frame. Using this
formula, we obtain some Bernstein theorems for $\Sigma$, which
generalize those results in \cite{han3} and \cite{han6} (see \S3
for details). Obviously, when $u_2$ and $u_3$ are constant,
$\Sigma$ is just the special Lagrangian graph in $C^n$. By
combining Wang's result in \cite{han3} and a Lewy transformation,
we establish a Bernstein Theorem for the special Lagrangian graph
under the condition: $Hess(u)\geq -CI$, where
$C<\frac{\sqrt{6}}{12}$. Note that our lower bound for $Hess(u)$
is independent of the dimension. Finally, we consider the minimal
Lagrangian graph given by three same potential functions, i.e.,
$\nabla u_s=\nabla u$ $(s=1,2,3)$. In this case, we also find a
suitable lewy transformation to prove a Bernstein result similar
to the above mentioned results.
\section{Preliminaries}
We first recall a formula derived in {\cite{han4}}, {\cite{han5}}.
Let $\Sigma$ be an oriented $n$-dimensional submanifold of
$R^{n+m}$ and $\Omega$ a parallel $n$ form on $R^{n+m}$. Around
any point $p\in\Sigma$, we choose any oriented orthonormal frames
$\{e_i\}_{i=1}^n$ for $T_p\Sigma$ and
$\{e_{\alpha}\}_{\alpha=n+1}^{n+m}$ for $N_p\Sigma$, the normal
bundle of $\Sigma$. The second fundamental form of $\Sigma$ is
denoted by $h_{\alpha ij}=<\nabla_{e_i}e_{j},e_{\alpha}>$. If we
assume $\Sigma$ has parallel mean curvature vector, then the
global function $*\Omega=\Omega(e_1,\cdots,e_n)$ satisfies (see
{\cite{han4}}, \cite{han5}):
\begin{eqnarray}
\triangle *\Omega+*\Omega(\sum_{\alpha lk}h_{\alpha,
l,k}^2)-2\sum_{\alpha,\beta, k}[\Omega_{\alpha\beta 3\cdots
n}h_{\alpha 1k}h_{\beta 2k}+\cdots+\Omega_{1\cdots
(n-2)\alpha\beta}h_{\alpha (n-1)k}h_{\beta nk}]=0 \label{a6}
\end{eqnarray}
\begin{eqnarray}
(*\Omega)_k=\sum_{\alpha}\Omega_{\alpha2\cdots n}h_{\alpha
1k}+\cdots+\Omega_{1\cdots(n-1)\alpha}h_{\alpha
nk}\label{a7}\end{eqnarray}
\begin{eqnarray*}
 1\leq i,j,k\cdots\leq n,\quad n+1\leq
\alpha,\beta\leq n+m
\end{eqnarray*}
where $\triangle$ is the Laplace operator of the induced metric on
$\Sigma$ and $\Omega_{\alpha\beta 3\cdots
n}=\Omega(e_{\alpha},e_{\beta},e_3,\cdots\\
\cdots,e_n)$, etc.

Let $g$ be the standard Euclidean metric on $H^n\cong R^{4n}$
which is K$\ddot{a}$hler with respect to the three natural complex
structures $I,J,K$ on $R^{4n}$.
 Set
$\omega_I=g(I,)$, $\omega_J=g(J,)$ and $\omega_K=g(K,)$. An
$n$-dimensional submanifold $f:\Sigma^n\longrightarrow R^{4n}$ is
called Lagrangian if it satisfies:
$$f^*\omega_I=f^*\omega_J=f^*\omega_K=0$$
Suppose $\Sigma$ is a graph defined by $f=(f_1,f_2,f_3)$, where
$f_s=(f_s^1,\cdots,f_s^n):R^n\longrightarrow R^n$ $s=1,2,3$ are
smooth maps. It is easy to see that $\Sigma$ is Lagrangian if and
only if $f_1,f_2,f_3$ satisfy:
\begin{eqnarray}
\left\{
\begin{array}{cccc}
-\frac{\partial f^j_1}{\partial x_i}+\frac{\partial
f_1^i}{\partial x_j}+\sum_{k=1}^n\{ \frac{\partial f_2^k}{\partial
x_i}\frac{\partial f_3^k}{\partial x_j}- \frac{\partial
f_2^k}{\partial x_j}\frac{\partial f_3^k}{\partial x_i}
\}&=0&\\
-\frac{\partial f^j_2}{\partial x_i}+\frac{\partial
f_2^i}{\partial x_j}+\sum_{k=1}^n\{ \frac{\partial f_3^k}{\partial
x_i}\frac{\partial f_1^k}{\partial x_j}- \frac{\partial
f_3^k}{\partial x_j}\frac{\partial f_1^k}{\partial x_i}
\}&=0&\\
-\frac{\partial f^j_3}{\partial x_i}+\frac{\partial
f_3^i}{\partial x_j}+\sum_{k=1}^n\{ \frac{\partial f_1^k}{\partial
x_i}\frac{\partial f_2^k}{\partial x_j}- \frac{\partial
f_1^k}{\partial x_j}\frac{\partial f_2^k}{\partial x_i}
\}&=0&\\
\end{array}\right.\label{a10}
\end{eqnarray}
where $i,j\in\{1,\cdots,n\}$. Obviously, if $f_s=\nabla u_s$ for
some smooth functions $u_s:R^n\longrightarrow R$ $(s=1,2,3)$, then
$\Sigma$ is Lagrangian if and only if $u_s$, $s=1,2,3$ satisfy:
\begin{eqnarray}Hess
(u_s)Hess (u_t)=Hess (u_t)Hess(u_s)\quad \text{for all}\quad
s,t\in\{1,2,3\}\label{a4}
\end{eqnarray}
From now on, we will consider minimal Lagrangian graph given by
the three potential functions $u_s$, $s=1,2,3$. Note that if
$u_2,u_3$ are constants, $\Sigma$ is just the special Lagrangian
graph considered in \cite{han2}, \cite{han3} and \cite{han6}.
\begin{exm}
The smallest interesting dimension is $n=2$. We would like to give
some examples of minimal Lagrangian surfaces in $H^2=R^8$.

First, let $u_s(x_1,x_2)$, $s=1,2,3$ be harmonic functions on
$R^2$. We see that $(u_s)_{x_1}-\sqrt{-1}(u_s)_{x_2}$ is a
holomorphic function of $z=x_1+\sqrt{-1}x_2$ for each $s$. It
follows that the graph $\Sigma=\{(x_1,x_2,\nabla u_1,\nabla
u_2,\nabla u_3):x=(x_1,x_2)\in R^2\}$ is a holomorphic curve in
$C^4=R^8$, and thus a minimal surfaces in $R^8$. In particular, if
$u $ is a harmonic function on $R^2$, we have minimal Lagrangian
graphs $\Sigma_1=\{(x,\nabla u,x,\nabla u):x\in R^2\}$ and
$\Sigma_2=\{(x,\nabla u,\nabla u,\nabla u):x\in R^2\}$.\label{a18}

\end{exm}

From Example \ref{a18}, we know that there exist many minimal
graphic Lagrangian submanifolds of $R^{4n}$.

\section{Main results for Minimal Lagrangian graphs }

Let $\Sigma=(x,\nabla u_1,\nabla u_2,\nabla u_3)$ be an $n$
dimensional minimal Lagrangian submanifold in $R^{4n}$. From the
previous section, we know that $\{u_s\}_{s=1,2,3}$ satisfy
(\ref{a4}) at each point $x\in R^n$. So we may diagonalize
$Hess(u_s)$, $s=1,2,3$ simultaneously at each point $x$ via the
singular decomposition, that is, there exist orthonormal bases
$\{a_i\}_{i=1,\cdots,n}$ for $R^n$ and
$\{a_{\alpha}\}_{\alpha=n+1,\cdots,4n}$ for $R^{3n}$ such that
$$Hess(u_s)a_i=\lambda_i^{(s)}a_{sn+i}$$
and
\begin{eqnarray*}
Ia_i=a_{n+i},\quad Ja_i=a_{2n+i},\quad Ka_i=a_{3n+i}
\end{eqnarray*}
for $i=1,\cdots,n$. Set
\begin{eqnarray*}
e_i&=&\frac{a_i+\lambda_i^{(1)}a_{n+i}+\lambda_i^{(2)}a_{2n+i}+\lambda_i^{(3)}a_{3n+i}}{A_i}\\
e_{n+i}&=&\frac{-\lambda_i^{(1)}a_i+a_{n+i}-\lambda_i^{(3)}a_{2n+i}+\lambda_i^{(2)}a_{3n+i}}{A_i}\\
e_{2n+i}&=&\frac{-\lambda_i^{(2)}a_i+\lambda_i^{(3)}a_{n+i}+a_{2n+i}-\lambda_i^{(1)}a_{3n+i}}{A_i}\\
e_{3n+i}&=&\frac{-\lambda_i^{(3)}a_i-\lambda_{i}^{(2)}a_{n+i}+\lambda_i^{(1)}a_{2n+i}+a_{3n+i}}{A_i}
\end{eqnarray*}
where
$A_i=\sqrt{1+(\lambda_i^{(1)})^2+(\lambda_i^{(2)})^2+(\lambda_i^{(3)})^2}$.
Note that, $e_{n+i}=Ie_i$, $e_{2n+i}=Je_i$ and $e_{3n+i}=Ke_i$ at
the corresponding point $p=(x,\nabla u_1(x),\nabla u_2(x),\nabla
u_3(x))$. Thus we have an orthonormal frame
$\{e_i\}_{i=1,\cdots,n}$ for $T_p\Sigma$ and
$\{e_{n+i},e_{2n+i},e_{3n+i}\}_{i=1,\cdots,n}$ for $N_p\Sigma$.
Define the second fundamental form of $\Sigma$ as follows:
$$h_{ijk}^{(s)}=<\widetilde{\nabla}_{e_i}e_j,e_{ns+k}>,\quad s=1,2,3$$
Since
$\widetilde{\nabla}I=\widetilde{\nabla}J=\widetilde{\nabla}K=0$
and $\Sigma$ is Lagrangian, we know that $h_{ijk}^{(s)}$ is
symmetric in $i,j,k$. Now take $\Omega=dx_1\wedge \cdots \wedge
dx_n$. It is not hard to see
$$*\Omega=\frac{1}{\sqrt{\Pi_{i=1}^n(1+\sum_{s=1}^3(\lambda_i^{(s)})^2)}}$$
By applying the formula (\ref{a6}) to the above quaternionic frame
$\{e_i,e_{n+i},e_{2n+i},e_{3n+i}\}_{i=1,\cdots, n}$, we get
\begin{prop}
Let $\Sigma=(x,\nabla u_1(x),\nabla u_2(x),\nabla u_3(x))$ be a
minimal graph in $R^{4n}$ and $\{\lambda_i^{(s)}\}$ be the
eigenvalues of $Hess(u_s)$, $s=1,2,3$. Then $*\Omega$ satisfies
\begin{eqnarray}
\triangle
*\Omega=-*\Omega\{\sum_{s=1}^3\sum_{ijk=1}^n(h_{ijk}^{(s)})^2-2\sum_{st=1}^3\sum_{k,i<j}
\lambda_i^{(s)}\lambda_j^{(t)}h_{iik}^{(s)}h_{jjk}^{(t)}+2\sum_{st=1}^3\sum_{k,i<j}\lambda_i^{s}\lambda_{j}^{(t)}
h_{ijk}^{(s)}h_{ijk}^{(t)}\}\label{a9}
\end{eqnarray}
where $\triangle$ is the Laplace operator of the induced metric on
$\Sigma$.
\end{prop}
Now we shall calculate
\begin{eqnarray}
\triangle
(ln*\Omega)=\frac{*\Omega\triangle(*\Omega)-|\nabla*\Omega|^2}{|*\Omega|^2}\label{a8}
\end{eqnarray}
By formula (\ref{a7}), the covariant derivative of $*\Omega$ is
\begin{eqnarray}
(*\Omega)_k=-*\Omega(\sum_{s=1}^3\sum_{i=1}^n\lambda_i^{(s)}h_{iik}^{s})
\end{eqnarray}
Plug this and equation (\ref{a9}) into equation (\ref{a8}) and we
obtain:
\begin{eqnarray}
\triangle
[ln(*\Omega)^{-1}]=\sum_{s=1}^3\sum_{ijk=1}^n(h_{ijk}^{(s)})^2+\sum_{st=1}^3\sum_{ijk=1}^n\lambda_i^{(s)}
\lambda_{j}^{(t)}h_{ijk}^{(s)}h_{ijk}^{(t)}
\end{eqnarray}
Set $\Lambda_i=(\lambda_i^{(1)},\lambda_i^{(2)},\lambda_i^{(3)})$
and $h_{ijk}=(h_{ijk}^{(1)},h_{ijk}^{(2)},h_{ijk}^{(3)})$. So
\begin{eqnarray}\triangle ln[*\Omega]^{(-1)}=\{\sum_{ijk=1}^{n}h_{ijk}(I+\Lambda_i^T
\Lambda_j)h_{ijk}^T\}\label{a16}
\end{eqnarray}
We may rewrite (\ref{a16}) as
\begin{eqnarray}\triangle ln[*\Omega]^{(-1)}=\frac{1}{3}\{\sum_{ijk=1}^{n}h_{ijk}(3I+\Lambda_i^T \Lambda_j
+\Lambda_j^T\Lambda_k+\Lambda_k^T\Lambda_i)h_{ijk}^T\}\label{a1}
\end{eqnarray}
Set
$S_{ij}=\frac{1}{2}(\Lambda_i^T\Lambda_j+\Lambda_j^T\Lambda_i)$.
We have the following:
\begin{thm}\label{a3}
Let $\Sigma=(x,\nabla u_1,\nabla u_2,\nabla u_3)$ be an
n-dimensional minimal Lagrangian submanifold of $R^{4n}$. If there
exist $\delta,K>0$ such that
$$|\lambda_i^{(s)}|\leq K,\quad\text{and}\quad
S_{ij}+S_{jk}+S_{ki}\geq (-3+\delta)I$$
for$i,j,k\in\{1,\cdots,n\}$, $s\in\{1,2,3\}$, then $\Sigma$ is an
affine plane.
\end{thm}
\begin{proof}
Set
\begin{eqnarray*}F_{ijk}(X)&=&X(3I+\Lambda_i^T \Lambda_j
+\Lambda_j^T\Lambda_k+\Lambda_k^T\Lambda_i)X^T\\
&=&X(3I+S_{ij}+S_{jk}+S_{ki})X^T\\
\end{eqnarray*}
for fixed $i,j,k$ and any $X\in R^3$. By the assumption we have
$$F_{ijk}(X)\geq\delta ||X||^2$$
From (\ref{a1}) we have
\begin{eqnarray*}\triangle ln[*\Omega]^{(-1)}
&=&\frac{1}{3}\sum_{ijk=1}^nF_{ijk}(h_{ijk})\\
&\geq&\frac{1}{3}\sum_{ijk=1}^n\delta||h_{ijk}||^2\\
&=&\frac{1}{3}\delta||A||^2
\end{eqnarray*}
where $A$ is the second fundamental form of $\Sigma$. Note that
$|\lambda_i^{(s)}|\leq K$ means $\Sigma$ is of bounded slope. So
we may perform blow down to get a minimal Lagrangian cone.
Obviously the minimal Lagrangian cone also satisfies the
assumption. By applying maximum principle we conclude that the
minimal cone is flat and then Allard regularity theorem implies
that $\Sigma$ is an affine plane.
\end{proof}

\begin{cor}\label{a22}
Let $\Sigma=(x,\nabla u_1,\nabla u_2,\nabla u_3)$ be a minimal
Lagrangian submanifold of $R^{4n}$. If there exist $\delta, K>0$
such that
$$|\lambda_i^{(s)}|\leq K,\quad\text{and}\quad
S_{ij}\geq (-\frac{3}{2}+\delta)I$$ where $i,j\in\{1,\cdots,n\}$,
$s\in\{1,2,3\}$, then $\Sigma$ is an affine plane.
\end{cor}
\begin{proof}
\begin{eqnarray*}F_{ijk}(X)&=&X(3I+\Lambda_i^T \Lambda_j
+\Lambda_j^T\Lambda_k+\Lambda_k^T\Lambda_i)X^T\\
&=&3||X||^2+XS_{ij}X^T+XS_{jk}X^T+XS_{jk}X^T\\
\end{eqnarray*}
It is easy to know that
\begin{eqnarray*}XS_{ij}X^T=(X,\Lambda_i)(\Lambda_j,X),\quad
XS_{jk}X^T=(X,\Lambda_j)(\Lambda_k,X),\quad
XS_{ki}X^T=(X,\Lambda_k)(\Lambda_i,X)
\end{eqnarray*}
Observe that one of $(X,\Lambda_i)(X,\Lambda_j),$
$(X,\Lambda_j)(X,\Lambda_k)$, $(X,\Lambda_k)(X,\Lambda_i)$ must be
nonnegative. From the assumption, we know
$$F_{ijk}(X)\geq 2\delta ||X||^2$$
i.e. $$S_{ij}+S_{jk}+S_{ki}\geq(-3+2\delta)I$$
 The conclusion follows immediately from
Theorem \ref{a3}.

\end{proof}
\begin{rem}When $u_1$ and $u_2$ are constants, we can recover Wang's result
 in \cite{han3} for special Lagrangian submanifolds.
\end{rem}

\begin{cor}\label{a11}
Let $\Sigma=(x,\nabla u_1,\nabla u_2,\nabla u_3)$ be a minimal
Lagrangian submnanifold of $R^{4n}$. If there exists a small
positive number $\delta$ such that
$$|\Lambda_i|\leq \sqrt{\frac{3}{2}-\delta}$$
for $i\in\{1,\cdots,n\}$, then $\Sigma$ is an affine plane.
\end{cor}
\begin{proof}
By the assumption $|\Lambda_i|\leq \sqrt{\frac{3}{2}-\delta}$ for
$i=1,\cdots,n$ and the Cauchy-Schwarz inequality, we have
$$XS_{ij}X^T=(X,\Lambda_i)(X,\Lambda_j)\geq (-\frac{3}{2}+\delta)||X||^2$$
for any fixed $i,j\in \{1,\cdots,n\}$ and $X\in R^3$. So the
symmetric matrix $S_{ij}$ satisfies
$$S_{ij}\geq (-\frac{3}{2}+\delta)I$$
Therefore the conclusion follows immediately from Corollary
\ref{a22}.
\end{proof}

\begin{rem}\label{h1}
If $u_2$ and $u_3$ are constants, then $\Sigma=(x,\nabla
u_1(x),0,0)$, which may be regarded as a minimal Lagrangian graph
$\Sigma=(x,\nabla u_1(x))$ in $C^{n}$. So the above Corollary
generalizes those results in \cite{han3} and \cite{han6}.
\end{rem}
In the following, we will consider two special kinds of minimal
Lagrangian graphs: $\Sigma=(x,\nabla u,0,0)$ or $\Sigma=(x,\nabla
u,\nabla u,\nabla u)$ in $H^n$. We have already pointed out that
the previous case is just the special Lagrangian graph.
\begin{thm}
Let $\Sigma=(x,\nabla u)$ be a minimal Lagrangian submanifold of
$C^n$. If there exists a positive constant $C<\frac{\sqrt{6}}{12}$
such that
$$Hess(u)\geq -CI$$
then $\Sigma$ is an affine plane.\label{y1}
\end{thm}
\begin{proof}
We identify $C$ with $R^2$ as follows:
$$C\ni x+\sqrt{-1}y\longleftrightarrow(x,y)\in R^2$$
For $a+\sqrt{-1}b\in SU(1)$, its real representation matrix on
$R^2$ is given by
\begin{eqnarray*}
A=\left(
\begin{array}{cccc}
a&&-b\\
b&&a
\end{array}\right)
\end{eqnarray*}
where $a^2+b^2=1$.

We consider the transformation $A^{(n)}=(A,\cdots,A)$ on
$\underbrace{C\times\cdots\times C}_n=R^{2n}$ defined by
\begin{eqnarray}
\left\{
\begin{array}{cccc}
\bar{x}&=&ax+by\\
\bar{y}&=&-bx+ay\\
\end{array}\right.
\end{eqnarray}
where $(x,y), (\bar{x},\bar{y})\in R^n\times R^n=C^n$. Set
$a=\frac{h}{\sqrt{1+h^2}}$ and $b=\frac{1}{\sqrt{1+h^2}}$, where
$h$ is a constant to be chosen. It follows that $\Sigma$ has a new
parametrization
\begin{eqnarray}
\left\{
\begin{array}{cccc}
\bar{x}&=&\frac{1}{\sqrt{1+h^2}}(h x+\nabla u)\\
\bar{y}&=&\frac{1}{\sqrt{1+h^2}}(-x+h\nabla u)\\
\end{array}\right.\label{h6}
\end{eqnarray}
Since $u+\frac{1}{2}C||x||^2$ is a convex function, we have
\begin{eqnarray}
|\bar{x}_1-\bar{x}_2|^2&=&\frac{1}{1+h^2}|hx_2+\nabla
u(x_2)-hx_1-\nabla u(x_1)|^2\nonumber\\
&=&\frac{1}{1+h^2}|(h-C)(x_2-x_1)+(\nabla u(x_2)+Cx_2)-(\nabla
u(x_1)+Cx_1)|^2\nonumber\\
&=&\frac{1}{1+h^2}\{(h-C)^2|x_2-x_1|^2+2(h-C)(x_2-x_1)[(\nabla
u(x_2)\nonumber\\
&+&Cx_2)-(\nabla u(x_1)+Cx_1)]+|(\nabla u(x_2)+Cx_2)-(\nabla
u(x_1)+Cx_1)|^2\}\nonumber\\
&\geq&\frac{1}{1+h^2}(h-C)^2|x_2-x_1|^2\label{h5}
\end{eqnarray}
Now we assume $h>C$. Then (\ref{h5}) implies that $\Sigma$ is
still a graph over the whole $\bar{x}$ space $R^n$. Further
$\Sigma$ is still a Lagrangian graph over $\bar{x}$, that means
$\Sigma$ has the representation $(\bar{x},\nabla
\bar{u}(\bar{x}))$ with a potential function $\bar{u}\in
C^{\infty}(R^n)$. we may derive from (\ref{h6}) that
$$Hess(\bar{u}(\bar{x}))=(hI+Hess(u(x)))^{-1}(-I+hHess(u(x)))$$
From $D^2u\geq-CI$, we see that
\begin{eqnarray}
-\frac{1+hC}{h-C}I\leq Hess(\bar{u})\leq hI\label{h7}
\end{eqnarray}
By solving $h=\frac{1+hC}{h-C}$, we get $h=C+\sqrt{C^2+1}$ which
obviously satisfies the previous assumption $h>C$. So (\ref{h7})
becomes
$$-(C+\sqrt{C^2+1})I\leq Hess(\bar{u}(\bar{x}))\leq(C+\sqrt{C^2+1})I$$
The condition $C\leq \frac{\sqrt{6}}{12}$ implies that any
eigenvalue $\lambda$ of $Hess(\bar{u})$ satisfies
$$|\lambda|\leq C+\sqrt{C^2+1}<\sqrt{\frac{3}{2}}$$
From Corollary \ref{a11}, we know that $\Sigma$ is an affine
plane.
\end{proof}
\begin{rem}
The lower bound for $Hess(u)$ in Theorem \ref{y1} is independent
of the dimension of Lagrangian graph $\Sigma=(x,\nabla u)$. This
improves Yuan's results in \cite{han6}.

\end{rem}
\begin{thm}
Let $\Sigma=(x,\nabla u,\nabla u,\nabla u)$ be a minimal
Lagrangian submanifold of $R^{4n}$. If there exists a positive
constant $C<\frac{\sqrt{2}}{12}$ such that
$$Hess(u)\geq -CI$$
then $\Sigma$ is an affine plane.
\end{thm}
\begin{proof}
We identify $H$ with $C^2$ as follows:
$$H\ni x+Iy+Jz+Kw=(x+Jz)+I(y+Jw)\longleftrightarrow (x+Jz,y+Jw)\in C^2$$
For a matrix $M=A+jB\in SU(2)$, its real representation on $R^4$
is given by
\begin{eqnarray*}
M\left(
\begin{array}{cccc}
x\\
y\\
z\\
w
\end{array}\right)=
 \left(
\begin{array}{cccc}
A&&-B\\
B&&A
\end{array}\right)
\left(
\begin{array}{cccc}
x\\
y\\
z\\
w
\end{array}\right)
\end{eqnarray*}
with $AA^T+BB^T=I_2$ and $AB^T=BA^T$. Obviously if $a^2+3b^2=1$,
then
\begin{eqnarray}
D=\left(
\begin{array}{ccccccccccccccccc}
a&&b&&b&&b\\
-b&&a&&b&&-b\\
-b&&-b&&a&&b\\
-b&&b&&-b&&a
\end{array}\right)\in SU(2)=Sp(1)
\end{eqnarray}
We consider the transformation $D^{(n)}=(D,\cdots,D)$ on
$\underbrace{H\times\cdots\times H}_n=H^n=R^{4n}$ defined by
\begin{eqnarray}
\left\{
\begin{array}{cccc}
\bar{x}&=&ax+by+bz+bw\\
\bar{y}&=&-bx+ay+bz-bw\\
\bar{z}&=&-bx-by+az+bw\\
\bar{w}&=&-bx+by-b z+aw\\
\end{array}\right.\label{a12}
\end{eqnarray}
where $(x,y,z,w),(\bar{x},\bar{y},\bar{z},\bar{w})\in R^n\times
R^n\times R^n \times R^n=R^{4n}$. Set $a=\frac{h}{\sqrt{1+h^2}}$
and $\sqrt{3}b=\frac{1}{\sqrt{1+h^2}}$, where $h$ is a constant to
be chosen. It follows that $\Sigma$ has a new parametrization
\begin{eqnarray}
\left\{
\begin{array}{cccc}
\bar{x}&=&\frac{1}{\sqrt{1+h^2}}(h x+\sqrt{3}\nabla u)\\
\bar{y}&=&\frac{1}{\sqrt{1+h^2}}(-\frac{1}{\sqrt{3}}x+h\nabla u)\\
\bar{z}&=&\frac{1}{\sqrt{1+h^2}}(-\frac{1}{\sqrt{3}}x+h\nabla u)\\
\bar{w}&=&\frac{1}{\sqrt{1+h^2}}(-\frac{1}{\sqrt{3}}x+h\nabla u)\\
\end{array}\right.\label{h3}
\end{eqnarray}
Since $u+\frac{1}{2}C||x||^2$ is a convex function, we have
\begin{eqnarray}
|\bar{x}_1-\bar{x}_2|^2&=&\frac{1}{1+h^2}|hx_2+\sqrt{3}\nabla
u(x_2)-hx_1-\sqrt{3}\nabla u(x_1)|^2 \nonumber\\
&=&\frac{1}{1+h^2}|(h-\sqrt{3}C)(x_2-x_1)+\sqrt{3}[(\nabla
u(x_2)+Cx_2)-(\nabla
u(x_1)+Cx_1)]|^2\nonumber\\
&=&\frac{1}{1+h^2}\{[(h-\sqrt{3}C)^2|x_2-x_1|^2]+2\sqrt{3}(h-\sqrt{3}C)(x_2-x_1)[(\nabla
u(x_2)\nonumber\\
&+&Cx_2)-(\nabla u(x_1)+Cx_1)]+3|(\nabla
u(x_2)+Cx_2)-(\nabla
u(x_1)+Cx_1)|^2\}\nonumber\\
&\geq&\frac{1}{1+h^2}(h-\sqrt{3}C)^2|x_2-x_1|^2\label{h2}
\end{eqnarray}
Now we assume $h>\sqrt{3}C$. Then (\ref{h2}) implies that $\Sigma$
is still a graph over the whole $\bar{x}-R^n$, that is $\Sigma$
has the representation
$(\bar{x},f_1(\bar{x}),f_2(\bar{x}),f_3(\bar{x}))$. Since $\Sigma$
is still minimal Lagrangian, we see from (\ref{a10}) that
$f=\nabla \bar{u}$, that is, $\Sigma=(\bar{x},\nabla
\bar{u},\nabla \bar{u},\nabla \bar{u})$ for some function
$\bar{u}\in C^{\infty}(R^n)$. We may derive from (\ref{h3}) that
$$Hess(\bar{u}(\bar{x}))=(hI+\sqrt{3}Hess(u(x)))^{-1}(-\frac{1}{\sqrt{3}}I+hHess(u(x)))$$
From $Hess(u)\geq-CI$, we see that
\begin{eqnarray}
-\frac{\frac{1}{\sqrt{3}}+hC}{h-\sqrt{3}C}I\leq Hess(\bar{u})\leq
\frac{h}{\sqrt{3}}I \label{h4}
\end{eqnarray}
By solving
$\frac{h}{\sqrt{3}}=\frac{\frac{1}{\sqrt{3}}+hC}{h-\sqrt{3}C}$, we
get $h=\sqrt{3}C+\sqrt{3C^2+1}$ which obviously satisfies the
previous assumption $h>\sqrt{3}C$. So (\ref{h4}) becomes
$$-\frac{\sqrt{3}C+\sqrt{3C^2+1}}{\sqrt{3}}I\leq Hess(\bar{u})\leq\frac{\sqrt{3}C+
\sqrt{3C^2+1}}{\sqrt{3}}I$$ The condition $C<\frac{\sqrt{2}}{12}$
implies that any eigenvalue $\lambda$ of $Hess(\bar{u}(\bar{x}))$
satisfies
$$|\lambda|\leq\frac{\sqrt{3}C+\sqrt{3C^2+1}}{\sqrt{3}}<\sqrt{\frac{1}{2}}$$
Then the singular values
$\Lambda_i=(\lambda_i,\lambda_i,\lambda_i)$ of
$\Sigma=(\bar{x},\nabla \bar{u}(\bar{x}),\nabla \bar{u}(\bar{x}),
\nabla \bar{u}(\bar{x}))$ satisfy:
$$|\Lambda_i|\leq \sqrt{3}C+\sqrt{3C^2+1}<\sqrt{\frac{3}{2}}\quad\text{for}\quad i\in\{1,\cdots,n\}$$
Hence, by Corollary \ref{a11}, we know that $\Sigma$ is an affine
plane.
\end{proof}
\begin{cor}
Let $\Sigma=(x,\nabla u,\nabla u,\nabla u)$ be a minimal
Lagrangian submanifold of $R^{4n}$. If $u$ is a smooth convex
function on $R^n$, then $\Sigma$ is an affine plane.
\end{cor}
{\bf{Acknowledgement}}: We would like to thank Professors Gu, C.H.
and Hu, H.S. for their valuable suggestions and constant
encouragement. We also thank Prof. Y.L. Xin for his helpful
comments.

Yuxin Dong, Yingbo Han and  Qingchun Ji\\

Institute of Mathematics

Fudan University

Shanghai, 200433, P. R. China

And

Key Laboratory of Mathematics

for Nonlinear Sciences (Fudan University),

Ministry of Education\\

Email address: dhj06630@163.com

\end{document}